\documentstyle{amsart}
\input{amssym.def}
\input{amssym.tex}

\newtheorem{theorem}{Theorem}

\newtheorem{corollary}[theorem]{Corollary}

\theoremstyle{definition}

\newtheorem{definition}[theorem]{Definition}

\newtheorem{example}[theorem]{Example}

\newtheorem{remark}[theorem]{Remark}

\def\ord{{\rm{ord}}}

\begin{document}
\bibliographystyle{plain}
\title[Dynamical zeta functions]
{Dynamical zeta functions for typical extensions of full shifts}
\author{T. Ward}
\address{School of Mathematics\\
University of East Anglia\\
Norwich NR4 7TJ\\
U.K.}
\email{t.ward@@uea.ac.uk}
\date{18 January 1999}

\begin{abstract}{We consider a family of
isometric extensions of the full shift on $p$ symbols
(for $p$ a prime) parametrized by
a probability space. Using Heath--Brown's work on the
Artin conjecture, it is shown that for all but two
primes $p$ the set of limit points of the
growth rate of periodic points is infinite almost
surely.
This shows in particular that the dynamical
zeta function is not algebraic almost surely.}
\end{abstract}
\maketitle

\section{Introduction}

The $S$-integer dynamical systems were introduced
in \cite{chothi-everest-ward-1997}: they are a
natural family of isometric extensions of
hyperbolic dynamical systems, parametrized by
rings of $S$-integers in $\Bbb A$-fields.
Their dynamical properties are governed by
arithmetic in algebraic number fields or
rational function fields depending on the
characteristic. A detailed description
is in \cite{chothi-everest-ward-1997},
along with some examples;
the ``random'' approach to their study
is outlined in \cite{ward-1997-uncountable-family}
and \cite{ward-1997-typical-behaviour}. Applications to
a certain class of cellular automata are described
in \cite{ward-cellular-automata-1998}.

Our purpose here is to extend a result
from \cite{ward-1997-typical-behaviour} concerning typical
behaviour for simple examples associated to
$\Bbb A$-fields of finite characteristic.

Let $K=\Bbb F_p(t)$, and following
Weil \cite[Chapter III]{weil-1974-number} let $P=\{\vert\cdot\vert_{v}\}$ be
the set of places (equivalence classes
of inequivalent multiplicative valuations) on
$K$.
The ``finite'' elements of $P$ are in one-to-one
correspondence with the irreducible
polynomials in $\Bbb F_p[t]$,
with $\vert f\vert_{v(t)}=
p^{-\ord_{v}(f)\cdot\deg(f)}$.
The ``infinite'' place is
determined by $\vert f(t)\vert_{\infty}=
\vert f(t^{-1})\vert_t$. With these
normalizations the product formula
\begin{equation}
\label{productformula}
\prod_{v\le\infty}\vert f\vert_{v}=1\mbox{ for all }
f\in K\backslash\{0\}
\end{equation}
holds.
Enumerate the countable set $P$ in some
order
\begin{equation}
\label{enumerationofsetofplaces}
P=\{\vert\cdot\vert_{v_{-1}=\infty},
\vert\cdot\vert_{v_0=t},
\vert\cdot\vert_{v_1},\dots\}.
\end{equation}
Denote by $\Omega$ the probability space
$\{0,1\}^{\Bbb N}$, equipped with the infinite
product measure $\mu_{\rho}=(\rho,1-\rho)^{\Bbb N}$ for
some $\rho\in[0,1]$.
Each point $\omega=\left(
\omega(k)\right)_{k\ge1}\in\Omega$ defines
a ring $R_{\omega}$ of
$S$-integers in $k$ defined by
\begin{equation}
\label{ringdeterminedbyomega}
R_{\omega}=\{f\in K:
\vert f\vert_{v_k}\le1\mbox{ for all }
k\ge1\mbox{ such that }\omega(k)=0\}.
\end{equation}
Notice that the infinite place and the place
corresponding to the irreducible polynomial $t$ are excluded from
the condition imposed in (\ref{ringdeterminedbyomega}).

\begin{example}\label{twoextremeexamplesofrings}
If $\rho=1$ then $\mu_{\rho}$-almost surely
$\omega(k)=0$ for all $k\ge 1$ so
$$
R_{\omega}=\{f\in K:\vert f\vert_{v_k}\le 1\mbox{ for all }k\ge 1\}
=\Bbb F_p[t^{\pm1}].
$$
At the opposite extreme, if $\rho=0$ then $\mu_{\rho}$-almost surely
$\omega(k)=1$ for all $k\ge 1$ so
$$
R_{\omega}=\Bbb F_p(t).
$$
\end{example}

For non-atomic measures $\mu_{\rho}$ (that is,
for $0<\rho<1$), the ring $R_{\omega}$
is a ``random'' ring, in which each irreducible
polynomial $v(t)$ is invertible with independent
probability $1-\rho$.

\begin{definition}\label{definitionofdynamicalsystem}
The dynamical system $\alpha_{\omega}:X_{\omega}\to
X_{\omega}$ is the automorphism $\alpha_{\omega}$ of
the compact abelian group $X_{\omega}$ dual to
the automorphism $f\mapsto tf$ of the ring
$R_{\omega}$.
\end{definition}

A basic formula from \cite{chothi-everest-ward-1997}
gives the number of periodic points in such a dynamical
system: the points of period $n$ under the
map $\alpha_{\omega}$ comprise the
set
$$
F_n(\alpha_{\omega})=\{x\in X_{\omega}:\alpha_{\omega}^n(x)=x\},
$$
and \cite[Lemma 5.2]{chothi-everest-ward-1997} shows that
\begin{equation}
\label{numberofperiodicpointsformula}
\left\vert F_n(\alpha_{\omega})\right\vert=\vert t^n-1\vert_{\infty}
\times\prod_{\omega(k)=1}\vert t^n-1\vert_{v_k}=
p^n\times\prod_{\omega(k)=1}\vert t^n-1\vert_{v_k}.
\end{equation}
Equivalently, this may be written
\begin{equation}
\label{numberofperiodicpointsformulalogarithmic}
\log_p\left\vert F_n(\alpha_{\omega})\right\vert=
n-\sum_{\omega(k)=1}\ord_{v_k}(t^n-1)\cdot\deg(v_k),
\end{equation}
which indicates how the valuations {\sl excluded} from
the condition in (\ref{ringdeterminedbyomega})
reduce the number of points of given period depending
on how the corresponding irreducible polynomials
divide $t^n-1$.

As a measure of the regularity of the periodic point behaviour
of such a dynamical system, we have the
{\sl dynamical zeta function}
\begin{equation}
\label{definitionofdynamicalzetafunction}
\zeta_{\alpha_{\omega}}(z)=
\exp\sum_{n=1}^{\infty}
\vert F_n(\alpha_{\omega})\vert\times\frac{z^n}{n},
\end{equation}
and the set ${\cal L}(\alpha_{\omega})$ of limit points of
the set $\left\{\frac{1}{n}
\log\vert F_n(\alpha_{\omega})\vert\right\}_{n\in\Bbb N}$.
It is clear from (\ref{numberofperiodicpointsformulalogarithmic})
that the zeta function converges in the disc
$\vert z\vert<\frac{1}{p}$ and that
${\cal L}$ is a subset of $[0,\log p]$.

\begin{example}\label{examplesofdynamicalsystems}
If $\rho=1$ then ($\mu_{\rho}$-almost
surely) $R_{\omega}=\Bbb F_p[t^{\pm1}]$,
so
$$
X_{\omega}=\{0,1\dots,p-1\}^{\Bbb Z},
$$
the two-sided shift space on $p$ symbols.
The endomorphism $\alpha_{\omega}$ is
the left shift
$$
(\alpha_{\omega}(x))_r=x_{r+1}
$$
for $x=(x_r)\in X_{\omega}$.
It is clear that there are $p^n$ points of
period $n$ under $\alpha_{\omega}$, which is
confirmed by equation (\ref{numberofperiodicpointsformula}).
The zeta function is rational,
$\zeta(z)=\frac{1}{1-pz}$, and ${\cal L}=\{\log p\}$ is
a singleton.

The opposite example has $\rho=0$, so
($\mu_{\rho}$-almost
surely) $R_{\omega}=\Bbb F_p(t)$. Here
the group $X_{\omega}$ is extremely complicated
(it is isomorphic to the quotient of the
adele ring $\Bbb F_p(t)_{\Bbb A}$ by the usual
discrete embedded copy of $\Bbb F_p(t)$ -- see
\cite[Chapter IV, Section 2]{weil-1974-number}),
and $\vert F_n(\alpha_{\omega})\vert=1$ for
all $n\ge 1$. Once again the zeta function
is rational, $\zeta(z)=\frac{1}{1-z}$, and ${\cal L}=
\{0\}$ is a singleton.
\end{example}

Our result is motivated by two things:
an example from \cite{chothi-everest-ward-1997}
and a theorem from \cite{ward-1997-typical-behaviour}. The example
corresponds to the (non-random) ring
$$
R=\{f\in K:
\vert f\vert_{v}\le 1\mbox{ for all }v\neq t^{-1},t-1\}=
\Bbb F_p[t][\frac{1}{t-1}].
$$

\begin{example}\label{examplefromcewpaper}
(cf. \cite[Example 8.5]{chothi-everest-ward-1997})
The endomorphism $\alpha:X\to X$ dual to
$f\mapsto tf$ on the ring
$\Bbb F_p[t][\frac{1}{t-1}]$
has
\begin{equation}
\label{examplefromcewpaperequation}
{\cal L}(\alpha)=
\left\{
\left(1-\textstyle\frac{1}{q}\right)\log p:
q\in\Bbb N\backslash p\Bbb N\right\}\cup\{\log p\}.
\end{equation}
\end{example}

\begin{theorem}\label{heathbrownfrometdspaper}
Assume that $0<\rho<1$. Then,
with the possible exception of two primes $p$,
${\cal L}(\alpha_{\omega})\supset\{0,\log p\}$
and the zeta function of $\alpha_{\omega}$ is
irrational for $\mu$-almost every $\omega\in\Omega$.
\end{theorem}

This follows from \cite[Theorem 3]{ward-1997-typical-behaviour}.
What we prove here is that the 
infinite collection of limit points seen in
Example \ref{examplefromcewpaper} is typical for the
random family of dynamical systems.

\begin{theorem}\label{mainresult}
Assume that $0<\rho<1$. Then,
with the possible exception of two primes $p$,
${\cal L}(\alpha_{\omega})$ is an infinite set
containing $0$, $\log p$, and a sequence
converging to $\log p$,
for $\mu_{\rho}$-almost every $\omega\in\Omega$.
\end{theorem}

An element $\ell\in{\cal L}(\alpha_{\omega})$ corresponds to
a singularity at $e^{-\ell}$ for the dynamical zeta function;
it follows that Theorem \ref{mainresult} forces the
zeta function to be non-algebraic.

\begin{corollary} Assume that $0<\rho<1$.
Then, with the possible exception of two
primes $p$, the dynamical zeta function
of $\alpha_{\omega}$ is $\mu_{\rho}$-almost
surely not an algebraic function.
\end{corollary}

\section{Proof of Theorem}

There are three ingredients to the proof of
Theorem \ref{mainresult}. The first is a deep
result due to Heath--Brown \cite{heath-brown-1986} on
the Artin conjecture:
with the possible exception of two primes $p$,
$p$ is a primitive root mod $q$ for infinitely
many primes $q$.

The second is a trivial consequence of the
Borel--Cantelli lemma from probability:
if $0<\rho<1$ and $(n_j)$ is an increasing
sequence of integers, then
$\omega(n_j)=1$ for infinitely many values of
$j$ and $\omega(n_j)=0$ for infinitely
many values of $j$, for $\mu_{\rho}$-almost every $\omega\in\Omega$.
More generally, if $(A_j)$ is an infinite sequence of
subsets of $\Bbb N$ with bounded cardinality and with
$n_j=\min(A_j)\to\infty$, then
$$
A_j\subset\{n\in\Bbb N:
\omega(n)=1\}
$$
for infinitely many $j$, and
$$
A_j\subset\{n\in\Bbb N:
\omega(n)=0\}
$$
for infinitely many $j$ for
$\mu_{\rho}$-almost every $\omega\in\Omega$.

The third is an elementary fact from Galois theory (see
\cite[Theorem 2.47]{lidl-niederreiter} for instance):
if $q$ is prime, then the polynomial
$t^{q-1}+t^{q-2}+\dots+1$ splits over $\Bbb F_p$
into $(q-1)/r$ irreducible factors, where $r$ is the least
positive integer for which $p^{r}\equiv 1$ mod $q$.
Writing
\begin{equation}
\label{cyclotomiclidldefinition}
\pi_{n}(t)=\prod_{{\rm gcd}(s,n)=1}(t-\xi^s)
\end{equation}
for the $n$th cyclotomic polynomial (where
$\xi$ is a primitive $n$th root of unity and $\mbox{gcd}(n,p)=1$),
$\pi_n$ factorizes over $\Bbb F_p[t]$
into $\phi(n)/d$ distinct
irreducibles of the same degree, where $d$ is the least positive integer
such that $p^d\equiv 1$ mod $n$.

It will be notationally convenient to pass to subsequences
without using additional suffixes: the sequence $(n_j)$
below is progressively thinned out as the proof proceeds.

Turning to the proof, we can find an infinite
sequence $(n_j)$ of primes greater than $p$ with the property that
$$
t^{n_j}-1=(t-1)\pi_{n_j}(t),
$$
where
\begin{equation}\label{definitionofpisubnj}
\pi_{n_j}(t)=t^{n_j-1}+t^{n_j-2}+\dots+1
\end{equation}
is an irreducible polynomial
in $\Bbb F_p[t]$.
By (\ref{definitionofpisubnj}),
$$
t^{n_j}\equiv 1\mbox{ in }
\Bbb F_p[t]/\langle\pi_{n_j}\rangle,
$$
a field of order $p^{n_j-1}$.
It follows that
$$
t^{qn_j}\equiv 1\mbox{ in }
\Bbb F_p[t]/\langle\pi_{n_j}\rangle
$$
for any fixed prime $q\neq p$, so
\begin{equation}\label{inequalityonorder}
\vert t^{qn_j}-1\vert_{\pi_{n_j}}\le
p^{-(n_j-1)}.
\end{equation}
The first step is to refine
(\ref{inequalityonorder}) by computing the
exact order of $\pi_{n_j}$ in
$t^{qn_j}-1$; it is greater than or equal to one
by (\ref{inequalityonorder}).
Assume that
\begin{equation}\label{doubleorderinexpression}
t^{qn_j}-1=(t-1)\pi_{n_j}^2(t)\cdot h(t)
\end{equation}
for some polynomial $h\in\Bbb F_p[t]$.

Recall that the {\sl order} of a
polynomial $g\in\Bbb F_p[t]$ with
$g(0)\neq0$ is
defined to be the least positive integer $e$
for which $g$ divides $t^e-1$ (cf.
\cite[Lemma 3.1, Definition 3.2]{lidl-niederreiter}).
By Theorem 3.11 {\it ibid},
we have
\begin{equation}\label{lidlresultonorderofpisquared}
\mbox{order}(\pi_{n_j}^2)=
c_j\cdot p^d,
\end{equation}
where $d$ is the least positive integer with
$p^d\ge 2$ (so $d=1$).
On the other hand, (\ref{doubleorderinexpression})
shows that
\begin{equation}\label{otherhand}
\mbox{order}(\pi_{n_j}^2)\mbox{ divides }
qn_j.
\end{equation}
We deduce that $c_jp\vert qn_j$, so
$p\vert n_j$ and hence
$p=n_j$, which is a contradiction.
It follows that (\ref{doubleorderinexpression})
is impossible, so we obtain a refinement of
(\ref{inequalityonorder}):
\begin{equation}\label{refinedsharpinequalityonorder}
\vert t^{qn_j}-1\vert_{\pi_{n_j}}=
p^{-(n_j-1)},
\end{equation}
or equivalently the order of $\pi_{n_j}$ in
$t^{qn_j}-1$ is exactly one for all $j$.

This allows the number of points of period
$qn_j$ under $\alpha_{\omega}$ to be
approximately calculated $\mu_{\rho}$-almost surely for
infinitely many values of $j$.

Without loss of generality each $n_j$
exceeds $q$; factorize $t^{qn_j}-1$ into cyclotomic
factors
\begin{equation}
\label{firstfactorisationstep}
t^{qn_j}-1=\prod_{d\vert qn_j}\pi_{d}(t)=
\pi_{n_j}(t)\cdot\pi_q(t)\cdot\pi_1(t)
\cdot\pi_{qn_j}(t).
\end{equation}
The last factor splits into $\frac{\phi(qn_j)}{d}=
\frac{(q-1)(n_j-1)}{d}$ irreducible factors, where
$d$ is the least positive integer with $p^d\equiv 1$
mod $qn_j$.
By construction, the order of $p$ mod $n_j$ is
$(n_j-1)$, so $d\ge(n_j-1)$. It follows that
$\pi_{qn_j}$ splits into no more than $(q-1)$ irreducible
polynomials, each of degree no smaller than $(n_j-1)$.
By Borel--Cantelli, we may assume therefore (by passing
to a subsequence in $j$) that
\begin{equation}
\label{firstcontrolthing}
k\in\{r\in\Bbb N:
\vert\pi_{qn_j}\vert_{v_r}\neq 1\mbox{ for some }j\}\implies
\omega(k)=0
\end{equation}
$\mu_{\rho}$-almost surely.

The second and third factors $\pi_1\pi_q$ in 
(\ref{firstfactorisationstep}) are
fixed as $j$ varies so we can ignore them: there is
a constant $A>0$ for which
\begin{equation}
\label{secondcontrolthing}
A\le\prod_{\omega(k)=1}\vert\pi_1\cdot\pi_q\vert_{v_k}
\le 1
\end{equation}
for all $\omega\in\Omega$.

The first factor in (\ref{firstfactorisationstep}) is
by construction itself irreducible, and by
(\ref{refinedsharpinequalityonorder}) does not appear in
any of the other three factors.
This fact, together with
(\ref{firstcontrolthing}) and
(\ref{secondcontrolthing}) gives by
(\ref{numberofperiodicpointsformula})
\begin{equation}
\label{randomperiodicpointcalculation}
p^{qn_j}p^{-(n_j-1)}A\le\vert F_{qn_j}(\alpha_{\omega})\vert\le
p^{qn_j}p^{-(n_j-1)}
\end{equation}
along some infinite sequence of $j$'s,
$\mu_{\rho}$-almost surely.
This may be written
\begin{equation}
\label{betterrandomstatement}
\vert F_{qn_j}(\alpha_{\omega})\vert=
p^{(q-1)n_j}\cdot B_j,
\end{equation}
where $pA\le B_j\le p.$
The growth rate of periodic points along this
sequence is then
\begin{eqnarray*}
\lim_{j\to\infty}\frac{1}{qn_j}\log
\vert F_{qn_j}(\alpha_{\omega})\vert&=&
\lim_{j\to\infty}\frac{1}{qn_j}\log p^{(q-1)n_j}
+\lim_{j\to\infty}\frac{1}{qn_j}\log B_j\\
&=&\left(1-\frac{1}{q}\right)\log p.
\end{eqnarray*}
That is, after excluding two
possible values of $p$, for any prime $q$ distinct from $p$
we have constructed a sequence of times
showing that $\left(1-\frac{1}{q}\right)\log p$ lies
in ${\cal L}(\alpha_{\omega})$ for almost every $\omega$.
Together with Theorem \ref{heathbrownfrometdspaper}, this
proves Theorem \ref{mainresult}.

\begin{remark}(1) The exact determination of the set
${\cal L}(\alpha_{\omega})$ for almost every $\omega$
in the general case is open. Methods from
\cite{ward-1997-typical-behaviour} suggest that
there is a single subset $A\subset[0,\log p]$ with the
property that ${\cal L}(\alpha_{\omega})=A$ almost
surely, but it is not clear what $A$ is, save that
with few exceptions it contains $0$ and $\log p$ and
the infinite sequence constructed above.
In the characteristic zero
case \cite{ward-1997-typical-behaviour}
again shows that there is a single set $A$
which gives the limit points almost surely,
but it is even less accessible: indeed, no single
element of $A$ is known (cf. \cite{ward-1997-uncountable-family},
\cite{ward-1997-typical-behaviour}).

\noindent(2) That the dynamical zeta function is
typically irrational (second part of
Theorem \ref{heathbrownfrometdspaper}) would follow at once if it were
known that $\omega\neq\omega^{\prime}\implies
\zeta_{\alpha_{\omega}}\neq\zeta_{\alpha_{\omega^{\prime}}}$
(since there are only countably many rational zeta
functions by \cite{bowen-lanford-1970}). This is
clear in simple zero-characteristic examples (see
\cite{ward-1997-uncountable-family}), but it is not
clear whether this implication holds in the current setting.

\noindent(3) The product measure $\mu_{\rho}$ on $\Omega$
is used above for simplicity, but it is clear from the
proofs that all that is needed is some kind of Borel--Cantelli
property. Thus, for instance, if $\mu$ is any probability
measure on $\Omega$ that is positive on open sets,
invariant under the left shift action $T:\Omega\to\Omega$
defined by $T(\omega)_{k}=\omega_{k+1}$, and ergodic for $T$,
then Theorem \ref{mainresult} holds with respect to $\mu$.

\end{remark}

\end{document}